\documentclass[11pt]{article}
\usepackage{amsmath,amsfonts,amssymb,amsthm}
\newtheorem{teo}{Theorem}[section]
\newtheorem{rem}[teo]{Remark}
\newtheorem{coro}[teo]{Corollary}

\newtheorem{lem}[teo]{Lemma}
\theoremstyle{definition}
\newtheorem{defi}{Definition}[section]

\newcommand{\prue}{\noindent \textbf{Proof.~~}}
\newcommand{\R}{\mathbb{R}} 
\newcommand{\N}{\mathbb{N}} 

\newcommand{\fin}{ \hfill $\square$}

\newcommand{\tint}{\displaystyle\int}
\newcommand{\tlim}{\displaystyle\lim}

\newcommand{\defrac}{\displaystyle\frac}

\title{Altering distance functions and fixed point theorems through rational expression}

\author{J.R. Morales$^a$ and E.M. Rojas$^b$\\\\
\small $^a$Departamento de Matem\'aticas,\\
\small Universidad de Los Andes, M\'erida-5101, Venezuela.\\ 
\small moralesj@ula.ve\\
 \small $^b$Departamento de Matem\'aticas,\\
 \small Pontificia Universidad Javeriana, Bogot\'a, Colombia\\
 \small edixon.rojas@javeriana.edu.co
}
\date{}
\begin{document}
\maketitle

\begin{abstract}

In this paper, using  altering distance functions we obtain a generalization of the results due to B.K. Das and S. Gupta \cite{DG} as well as the results given by B. Samet and H. Yazid \cite{SY}. Moreover, we study the so-called  property $P$ for the contraction mappings considered in this article.
\end{abstract}

\section{Introduction and preliminary facts}

In 1984, M.S. Khan, M. Swalech and S. Sessa \cite{KSS}  expanded the research of the metric fixed point theory to a new category 
by introducing a control function which they called an \textit{altering distance function}.

\begin{defi}[\cite{KSS}]

A function $\psi: \R_{+}\longrightarrow \R_{+}:=[0,+\infty)$ is called an altering distance function if the following properties are satisfied:
\begin{enumerate}
\item[$(\Psi_1)$] $\psi(t)=0\Leftrightarrow t=0$.
\item[$(\Psi_2)$] $\psi$ is monotonically non-decreasing.
\item[$(\Psi_3)$] $\psi$ is continuous.
\end{enumerate}
By $\Psi$ we denote the set of all altering distance functions.
\end{defi}

Using those control functions the authors extend the  Banach Contraction Principle by  taking $\psi=Id,$ (the identity mapping), 
in the inequality contraction \eqref{eq1.1} of the following theorem.

\begin{teo}[\cite{KSS}]\label{teo1.2}
Let $(M,d)$ be a complete metric space, let $\psi\in \Psi$ and let $S: M\longrightarrow M$ be a mapping which satisfies the following inequality
\begin{equation}\label{eq1.1}
\psi[d(Sx, Sy)]\leq a\psi[d(x,y)]
\end{equation} for all $x,y\in M$ and for some $0<a<1$. Then. $S$ has a unique fixed point $z_0\in M$ and moreover for each $x\in M$, $\tlim_{n\rightarrow \infty}S^nx=z_0.$
\end{teo}

Fixed point theorems involving the notion of altering distance functions has been widely studied, see for instance \cite{BS,PM,SY} and references therein.

On the other hand, in 1975, B.K. Das and S. Gupta \cite{DG} proved the following result.

\begin{teo}\label{teo1.3}
Let $(M,d)$ be a metric space and let $S: M\longrightarrow M$ be a given mapping such that,
\begin{enumerate}
\item[(i)]
\begin{equation}\label{eq1.2}
d(Sx,Sy)\leq a d(x,y)+bm(x,y)
\end{equation}
 for all $x,y\in M,\,\, a>0,\,\, b>0,\,\, a+b<1$  where
\begin{equation}\label{eq1.3}
m(x,y)=d(y,Sy)\defrac{1+d(x,Sx)}{1+d(x,y)}
\end{equation} for all $x,y\in M$. 

\item[(ii)] For some $x_0\in M$, the sequence of iterates $(S^nx_0)$ has a subsequence $(S^{n_k}x_0)$ with
$\tlim_{k\rightarrow \infty}S^{n_k}x_0=z_0.$ Then $z_0$ is the unique fixed point of $S$.
\end{enumerate}
\end{teo}

It is important to indicate that B.K. Das and S. Gupta \cite{DG} did not hypothesize that $M$ is a complete metric space but they used this fact in their proof.

Finally, in this paper we will study the property introduced by  G.S. Jeong and B.E. Rhoades in \cite{JR} which they called \textit{the property} $P$ \textit{in metric spaces}:

\begin{defi}
Let $S$ be a self mapping of a metric space $(M,d)$ with a nonempty fixed point set $F(S)$. Then $S$ is said to satisfy
the property $P$ if $F(S)=F(S^n)$ for each $n\in \N$.
\end{defi}

An interesting fact about mappings satisfying the property $P$ is that they have no nontrivial periodic points. For more information about this property see e.g., \cite{JR,JeR,RA}.
The following lemma given by G.U. Babu and P.P. Sailaja \cite{BS} will be used in the sequel in order to prove our main results.

\begin{lem}\label{lem1.5}
Let $(M,d)$ be a metric space. Let $(x_n)$ be a sequence in $M$ such that
\begin{equation}\label{eq1.4}
\tlim_{n\rightarrow \infty} d(x_n, x_{n+1})=0.
\end{equation} 
If $(x_n)$ is not a Cauchy sequence in $M$, then there exist an $\varepsilon_0>0$ and sequences of integers positive $(m(k))$ and $(n(k))$ with
\begin{equation*}
m(k)>n(k)>k
\end{equation*}
 such that,
 \begin{equation*}
d(x_{m(k)},\, x_{n(k)})\geq \varepsilon_0,\quad d(x_{m(k)-1},\, x_{n(k)})<\varepsilon_0
 \end{equation*}
 and
\begin{enumerate}
\item[(i)] $\tlim_{k\rightarrow \infty}d(x_{m(k)-1},\, x_{n(k)+1})=\varepsilon_0$,
\item[(ii)] $\tlim_{k\rightarrow \infty}d(x_{m(k)},\, x_{n(k)})=\varepsilon_0$,
\item[(iii)] $\tlim_{k\rightarrow \infty}d(x_{m(k)-1},\, x_{n(k)})=\varepsilon_0$.
\end{enumerate}
\end{lem}

\begin{rem}\label{ref1.6}
From Lemma \ref{lem1.5} is easy to get
\begin{equation*}
\tlim_{k\rightarrow \infty} d(x_{m(k)+1},\, x_{n(k)+1})=\varepsilon_0.
\end{equation*}
\end{rem}
In this paper, we consider the altering distance functions to generalize the results given in \cite{DG} and also we will obtain an extension of Theorem 2 given in \cite{SY}.
Moreover, we study the property $P$ for the contraction mappings considered in this work.

\section{Fixed point theorems}
\setcounter{equation}{0}

This section is devoted to generalize the Theorem \ref{teo1.3}, as well as generalize Theorem 2 of  \cite{SY}.

\begin{teo}\label{teo2.1}
Let $(M,d)$ be a complete metric space, let $\psi\in \Psi$ and let $S: M\longrightarrow M$ be a mapping which satisfies the following condition:
\begin{equation}\label{eq2.1}
\psi[d(Sx,Sy)]\leq a\psi[d(x,y)]+b\psi(m(x,y))
\end{equation}
 for all $x,y\in M,\, a>0,\, b>0,\, a+b<1$ and $m(x,y)$ is given by \eqref{eq1.3}. Then $S$ has a unique fixed point $z_0\in M$, and moreover for each $x\in M$ $\tlim_{n\rightarrow \infty} S^nx=z_0$.
\end{teo}

\prue
Let $x\in M$ be an arbitrary point and let $(x_n)$ be a sequence defined as follows: $x_{n+1}=Sx_n=S^{n+1}x,$ for each $n\geq 1$.
Now,
\begin{equation}\label{eq2.2}
\begin{array}{rcl}
\psi[d(x_n,\, x_{n+1})] &=& \psi[d(Sx_{n-1},\, Sx_n)]\\ 
&\leq& a\psi[d(x_{n-1},\, x_n)]+
b\psi(m(x_{n-1},\, x_n))
\end{array}
\end{equation} 
using \eqref{eq1.3},
\begin{equation*}
m(x_{n-1},\, x_n)=d(x_n,\, x_{n+1})\defrac{1+d(x_{n-1},\, x_n)}{1+d(x_{n-1},\, x_n)}=d(x_n,\, x_{n+1})
\end{equation*}
substituting it into \eqref{eq2.2}, one obtains
\begin{equation*}
\psi[d(x_n,\, x_{n+1})]\leq a\psi(d(x_{n-1},\, x_n))+b\psi(d(x_n,\, x_{n+1}))
\end{equation*} 
it follows that,
\begin{equation}
\begin{array}{rcl}
\psi[d(x_n,\, x_{n+1})] &\leq& \defrac{a}{1-b}\psi(d(x_{n-1},\, x_n))\\  
&\leq& \left(\defrac{a}{1-b}\right)^2 \psi(d(x_{n-2},\, x_{n-1}))\leq \ldots\\
&\leq& \left(\defrac{a}{1-b}\right)^n\psi(d(x_0,\, x_1)). \label{eq2.3}
\end{array}
\end{equation} 
Since $\frac{a}{1-b}\in (0,1)$, from \eqref{eq2.3} we obtain
\begin{equation*}
\tlim_{n\rightarrow \infty} \psi[d(x_n,\, x_{n+1})]=0.
\end{equation*}
From the fact that $\psi\in \Psi$, we have
\begin{equation}\label{eq2.4}
\tlim_{n\rightarrow \infty} d(x_n,\, x_{n+1})=0.
\end{equation}
 Now, we will show that $(x_n)$ is a Cauchy sequence in $M.$ Suppose that $(x_n)$ is not a Cauchy sequence, which means that there is a 
 constant $\varepsilon_0>0$ such that for each positive integer $k,$ there are positive integers $m(k)$ and $n(k)$ with $m(k)>n(k)>k$ such that
\begin{equation*}
d(x_{m(k)},\, x_{n(k)})\geq \varepsilon_0,\,\, d(x_{m(k)-1},\, x_{n(k)})<\varepsilon_0.
\end{equation*}
From Lemma \ref{lem1.5} and Remark \ref{ref1.6} we obtain

\begin{equation}\label{eq2.5}
\tlim_{k\rightarrow \infty} d(x_{m(k)},\, x_{n(k)})=\varepsilon_0
\end{equation} 
and
\begin{equation}\label{eq2.6}
\tlim_{k\rightarrow \infty} d(x_{m(k)+1},\, x_{n(k)+1})=\varepsilon_0.
\end{equation} 
For $x=x_{m(k)}$ and $y=y_{n(k)}$ from \eqref{eq2.1} we have,
\begin{equation*}
\begin{array}{rcl}
\psi[d(x_{m(k)+1},\, x_{n(k)+1})] &=& \psi[d(Sx_{m(k)},\, x_{n(k)})]\leq a\psi[d(x_{m(k)},\, x_{n(k)})]\\
&&+ b\psi\left[d(x_{n(k)},\, x_{n(k)+1})\defrac{1+d(x_{m(k)},\,x_{m(k)+1})}{1+d(x_{m(k)},\, x_{n(k)})}\right]
\end{array}
\end{equation*}
 using \eqref{eq2.4}, \eqref{eq2.5} and \eqref{eq2.6} we obtain
\begin{equation*}
\begin{array}{rcl}
\psi(\varepsilon) &=& \tlim_{k\rightarrow \infty} \psi[d(x_{m(k)+1},\, x_{n(k)+1})]\\ 
 &\leq& a\tlim_{k\rightarrow \infty} \psi[d(x_{m(k)},\, x_{n(k)})]\\ 
 &\leq& a\psi(\varepsilon),
\end{array}
\end{equation*} 
since $a\in (0,1)$, we get a contradiction. Thus $(x_n)$ is a Cauchy sequence in the complete metric space $M$, thus there exists $z_0\in M$ such that
\begin{equation*}
\tlim_{n\rightarrow \infty} x_n=z_0.
\end{equation*}
Setting $x=x_n$ and $y=z_0$ in \eqref{eq2.1} we have
\begin{equation*}
\begin{array}{rcl}
\psi[d(x_{n+1},\, Sz_0)] &=& \psi[d(Sx_n,\, Sz_0)]\\ &\leq& a\psi[d(x_{n},\, z_0)]+ b\psi\left[d(z_0,\, Sz_0)\defrac{1+d(x_n,\,Sx_n)}{1+d(x_n,\, z_0)}\right].
\end{array}
\end{equation*}
 Therefore,
\begin{equation*}
\tlim_{n\rightarrow \infty} \psi[d(x_{n+1},\, Sz_0)]\leq b\psi[d(z_0,Sz_0)]
\end{equation*}
i.e.,
\begin{equation*}
\psi[d(z_0,\, Sz_0)]\leq b\psi[d(z_0,Sz_0)]
\end{equation*}
 since $b\in (0,1)$, then $\psi[d(z_0,Sz_0)]=0$ which implies that $d(z_0,Sz_0)=0$ thus $z_0=Sz_0$.

Now we are going to establish the uniqueness of the fixed point. 
Let $y_0, z_0$ be two fixed points of $S$ such that $y_0\neq z_0$. Putting $x=y_0$ and $y=z_0$ in \eqref{eq2.1} we get
\begin{eqnarray*}
\psi[d(Sz_0,\, Sy_0)]&\leq& a\psi[d(z_0,\, y_0)]+ b\psi\left[d(y_0,Sy_0)\defrac{1+d(z_0,Sz_0)}{1+d(z_0,\, y_0)}\right]\\
&=&a\psi[d(z_0,y_0)],
\end{eqnarray*}
 which implies that $\psi[d(z_0,\, y_0)]=0$, so $d(z_0,\, y_0)=0$. Thus $z_0=y_0$. \fin

\begin{coro}[\cite{SY}, Theorem 2]

Let $(M,d)$ be a complete metric space and let $T: M\longrightarrow M$ be a mapping. We assume that for each $x,y\in M,$

\begin{equation}\label{eq2.7}
\tint\limits_{0}^{\scriptscriptstyle d(Sx,Sy)} \varphi(t)dt\leq \tint\limits_{0}^{d(x,y)} \varphi(t)dt+b\tint\limits_{0}^{\scriptscriptstyle d(y,Sy)\frac{1+d(x,Sx)}{1+d(x,y)}} \varphi(t)dt
\end{equation} 
where $0<a+b<1$ and $\varphi: \R_{+}\longrightarrow \R_{+}$ is a Lebesgue integrable mapping which is summable on each compact subset of $[0,+\infty)$, non negative and such that $\int\limits_{0}^{\varepsilon}\varphi(t)dt>0$ for all $\varepsilon>0$. 
Then $S$ admits a unique fixed point $z_0\in M$ such that for each $x\in M$
\begin{equation*}
\tlim_{n\rightarrow \infty} S^nx=z_0.
\end{equation*}
\end{coro}
\prue Let $\varphi: \R_{+}\longrightarrow \R_{+}$ be as in the hypothesis,
 we define $\psi_0(t)=\tint_{0}^{t}\varphi(t)dt,\,\, t\in \R_{+}.$ It is clear that $\psi_0(0)=0$. $\psi_0$ is monotonically non decreasing
  and by hypothesis $\psi_0$ is absolutely continuous, hence $\psi_0$ is continuous. Therefore, $\psi_0\in \Psi$. So \eqref{eq2.1} becomes
\begin{equation*}
\psi_0(d(Sx,Sy))\leq a\psi_0(d(x,y))+b\psi_0\left[d(y,Sy)\defrac{1+d(x,Sx)}{1+d(x,y)}\right].
\end{equation*}
Hence from Theorem \ref{teo2.1} there exists a unique fixed point $z_0\in M$ such that for each $x\in M$, $\tlim_{n\rightarrow \infty} S^nx=z_0$. \fin

\begin{rem} ~\newline
\begin{enumerate}

\item[(1)] If we take $b=0$, then \eqref{eq2.1} reduces to \eqref{eq1.2}, thus the Theorem \ref{teo1.2} is a corollary
of Theorem \ref{teo2.1}.

\item[(2)] If we take $\psi=Id$ in \eqref{eq2.1}, then we obtain \eqref{eq1.2}. Therefore
the Theorem \ref{teo2.1} is a generalization of Theorem \ref{teo1.3}.
\end{enumerate}
\end{rem}

\section{The property $P$.}

In this section we are going to prove that the mappings satisfying the contractive conditions\eqref{eq1.1}, \eqref{eq1.2}, \eqref{eq2.1} and \eqref{eq2.7} fulfill the property $P$.

\begin{teo}

Let $(M,d)$ be a complete metric space, let $\psi\in \Psi$ and let $S: M\longrightarrow M$ be a
 mapping which satisfies the following inequality:
\begin{equation*}
\psi[d(Sx,Sy)]\leq a\psi[d(x,y)]
\end{equation*}
for all $x,y\in M$ and for some $0<a<1$. Then $F_{S}\neq \emptyset$ and $S$ has the property $P$.
\end{teo}

\prue From Theorem \ref{teo1.2}, $S$ has a fixed point. Therefore $F_{S^n}\neq \emptyset$ for each $n\in \N$.
Fix $n>1$ and we assume that $z\in F_{S^n}$ we want to show that $z\in F_{S}.$ Suppose that $z\neq Sz$, using \eqref{eq1.1}
\begin{equation*}
\begin{array}{rcl}
\psi[d(z,\, Sz)] &=& \psi[d(S^nz,\, S^{n+1}z)]\leq a\psi[d(S^{n-1}z,\, S^nz)]\\  &\leq& \ldots \leq a^n\psi[d(z,\, Sz)],
\end{array}
\end{equation*}
 since $a\in (0,1),$ $\tlim_{n\rightarrow \infty} \psi[d(z,\, Sz)]=0$. From the fact that $\psi\in \Psi,$ we get $ z=Sz$ which is a 
 contraction. Therefore $z\in F_{S}$ i.e., $S$ has the property $P$.\fin

\begin{teo}
Let $(M,d)$ be a complete metric space and let $S: M\longrightarrow M$ be a mapping which satisfies the contractive condition \eqref{eq1.2}, then $F_{S}\neq \emptyset$ and $S$ has the property $P$.
\end{teo}

\prue From Theorem \ref{teo1.3}, $F_{S}\neq \emptyset.$ Therefore $F_{S^n}\neq \emptyset$ for each $n\in \N$. Fix $n>1$ and we assume that 
$z\in F_{S^n}$. We want to show that $z\in F_{S}$. Suppose that $z\neq Sz$. Using \eqref{eq1.2},
\begin{equation*}
\begin{array}{rcl}
d(z,\, Sz) &=& d(S^nz,\, S^{n+1}z)\leq ad(S^{n-1}z,\, S^nz)\\
&&+bd(S^nz,\, S^{n+1}z)\defrac{1+d(S^{n-1}z,\, S^nz)}{1+d(S^{n-1}z,\, S^nz)}\\ 
 &=& ad(S^{n-1}z,\, S^nz)+bd(S^nz,\, S^{n+1}z).
\end{array}
\end{equation*}
Therefore
\begin{equation*}
d(z,Sz) = d(S^nz,\, S^{n+1}z)\leq \defrac{a}{1-b}d(S^{n-1}z,\, S^nz)\leq \ldots\leq \left(\defrac{a}{1-b}\right)^nd(z,Sz),
\end{equation*}
which is a contradiction. Consequently, $z\in F_{S}$ and $S$ has the property $P$. \fin

\begin{teo}

Let $(M,d)$ be a complete metric space, let $\psi\in \Psi$ and let $S: M\longrightarrow M$ be a mapping which satisfies the condition \eqref{eq2.1}. Then $F_{S}\neq \emptyset$ and $S$ has a property $P.$
\end{teo}

\prue From Theorem \ref{teo2.1}, $S$ has a fixed point. Therefore $F_{S^n}\neq \emptyset$ for each $n\in \N$. 
Fix $n>1$ and assume that $z\in F_{S^ n}$ we wish to show that $z\in F_{S}$. Suppose that $z\neq Sz$, using \eqref{eq2.1},
\begin{equation*}
\begin{array}{rcl}
\psi[d(z,Sz)] &=& \psi[d(S^nz,\, S^{n+1}z)]\leq a\psi[d(S^{n-1}z,\, S^nz)]\\
&&+b\psi\left[d(S^nz,\, S^{n+1}z)\defrac{1+d(S^{n-1}z,\, S^nz)}{1+d(S^{n-1}z,\, S^nz)}\right]\\ 
&\leq& a\psi[d(S^{n-1}z,\, S^nz)]+b\psi[d(S^nz,\, S^{n+1}z)]
\end{array}
\end{equation*}
hence
\begin{equation*}
\begin{array}{rcl}
(1-b)\psi[d(S^nz,\, S^{n+1}z)] &\leq& a\psi[d(S^{n-1}z,\, S^nz)]\\ 
 \psi[d(z,\, Sz)]=\psi(d(S^nz,\, S^{n+1}z)) &\leq& \defrac{a}{1-b} \psi[d(S^{n-1}z,\, S^nz)]\\ 
  &\leq& \left(\defrac{a}{1-b}\right)^n\psi(d(z,\, Sz)),
\end{array}
\end{equation*}
thus
\begin{equation*}
\psi[d(z,\, Sz)]\leq \left(\defrac{a}{1-b}\right)^n\psi(d(z,\, Sz))
\end{equation*}
 which is a contradiction, therefore $\psi[d(z,\, Sz)]=0$, since $\psi\in\Psi$, we conclude that $d(z, Sz)=0,$ thus $z\in F_{S}$ and $S$ has the property $P$.\fin

\end{document}